  \def\version{September 23, 2003}

\newif\ifpdf
\ifx\pdfoutput\undefined
\pdffalse 
\else
\pdfoutput=1 
\pdftrue \fi

\newif\iffinal
\finalfalse 
\finaltrue 

\documentclass[reqno,twoside,10pt]{amsart}
\iffinal\else\usepackage[notref,notcite]{showkeys}\fi
\usepackage{amsmath}
\usepackage{amsfonts}
\usepackage{amssymb}
\usepackage{verbatim}
\usepackage{enumerate}
\usepackage{color}
\usepackage[colorlinks,citecolor=blue,urlcolor=blue]{hyperref}
\IfFileExists{myowntimes.sty}{\usepackage{myowntimes}}{\usepackage{times}\usepackage{mathrsfs}}

\ifpdf
\usepackage[pdftex]{graphicx}
\else
\usepackage{epsfig}
\fi

\ifpdf \DeclareGraphicsExtensions{.pdf, .jpg} \else
\DeclareGraphicsExtensions{.eps, .jpg} \fi

\DeclareFontFamily{OT1}{eusb}{} \DeclareFontShape{OT1}{eusb}{m}{n}
{<5> <6> <7> <8> <9> <10> <11> <12> <14.4> eusb10}{}
\DeclareMathAlphabet{\eusb}{OT1}{eusb}{m}{n}

\DeclareFontFamily{OT1}{eusm}{} \DeclareFontShape{OT1}{eusm}{m}{n}
{<5> <6> <7> <8> <9> <10> <11> <12> <14.4> eusm10}{}
\DeclareMathAlphabet{\eusm}{OT1}{eusm}{m}{n}

\DeclareFontFamily{OT1}{eufm}{} \DeclareFontShape{OT1}{eufm}{m}{n}
{<5> <6> <7> <8> <9> <10> <11> <12> <14.4> eufm10}{}
\DeclareMathAlphabet{\mathfrak}{OT1}{eufm}{m}{n}

\DeclareFontFamily{OT1}{fraktura}{}
\DeclareFontShape{OT1}{fraktura}{m}{n} {<5> <6> <7> <8> <9> <10>
<11> <12> <13> <14.4> [1.1] eufm10}{}
\DeclareMathAlphabet{\fraktura}{OT1}{fraktura}{m}{n}

\DeclareFontFamily{OT1}{cmfi}{} \DeclareFontShape{OT1}{cmfi}{m}{n}
{<5> <6> <7> <8> <9> <10> <11> <12> <13> <14.4> [0.9] cmfi10}{}
\DeclareMathAlphabet{\cmfi}{OT1}{cmfi}{b}{n}

\newtheoremstyle{thm}{1.5ex}{1.5ex}{\itshape\rmfamily}{}
{\bfseries\rmfamily}{}{2ex}{}

\newtheoremstyle{rem}{1.3ex}{1.3ex}{\rmfamily}{}
{\itshape}
{} {1.5ex}{}

\newenvironment{proofsect}[1]
{\vskip0.1cm\noindent{\rmfamily\itshape #1.}}{\qed\vglue0.3cm}

\theoremstyle{thm}
\newtheorem{theorem}{Theorem}[section]
\newtheorem{lemma}[theorem]{Lemma}

\newtheorem*{Main Theorem}{Main Theorem}

\newtheorem{definition}{Definition}

\theoremstyle{rem}

\numberwithin{equation}{section}

\renewcommand{\section}{\secdef\sct\sect}
\newcommand{\sct}[2][default]{\refstepcounter{section}
\addcontentsline{toc}{section} {{\tocsection
{}{\thesection}{\!\!\!\!#1\dotfill}}{}} \vspace{0.7cm}
\centerline{ 
\scshape\arabic{section}.\ #1} \nopagebreak \vspace{0.2cm}}
\newcommand{\sect}[1]{
\vspace{0.4cm} \centerline{\large\scshape\rmfamily #1}
\vspace{0.2cm}}

\renewcommand{\subsection}{\secdef\subsct\sbsect}
\newcommand{\subsct}[2][default]{\refstepcounter{subsection}
\addcontentsline{toc}{subsection}
{{\tocsection{\!\!}{\hspace{1.2em}\thesubsection}{\!\!\!\!#1\dotfill}}{}}
\nopagebreak \vspace{0.45\baselineskip} {\flushleft\bf
\arabic{section}.\arabic{subsection}~\bf #1.~}
\\*[3mm]\noindent
\nopagebreak}
\newcommand{\sbsect}[1]{\vspace{0.1cm}\noindent
\textbf{#1.~}\vspace{0.1cm}}

\renewcommand{\subsubsection}{%
\secdef \subsubsect\sbsbsect}
\newcommand{\subsubsect}[2][default]{%
\refstepcounter{subsubsection}
\addcontentsline{toc}{subsubsection}{{\tocsection{\!\!}
{\hspace{3.05em}\thesubsubsection}{\!\!\!\!#1\dotfill}}{}}
\nopagebreak \vspace{0.15\baselineskip} \nopagebreak
{\flushleft\rmfamily
\itshape\arabic{section}.\arabic{subsection}.\arabic{subsubsection}
\ \rmfamily #1\/.}\ }
\newcommand{\sbsbsect}[1]{\vspace{0.1cm}\noindent
\rmfamily \itshape
\arabic{section}.\arabic{subsection}.\arabic{subsubsection} \
\sffamily #1\/.\ }

\iffinal
\newcommand{\printversion}{}
\else
\newcommand{\printversion}{, \version}
\fi

\begin{document}

\title[THE TIGHTNESS OF THE KESTEN-STIGUM RECONSTRUCTION BOUND FOR A
SYMMETRIC MODEL WITH MULTIPLE MUTATIONS]{%
\Large THE TIGHTNESS OF THE KESTEN-STIGUM RECONSTRUCTION BOUND FOR A
SYMMETRIC MODEL WITH MULTIPLE MUTATIONS}
\author[ Wenjian~Liu, S. RAO JAMMALAMADAKA and Ning Ning\printversion]{Wenjian~Liu$^*$, S. RAO JAMMALAMADAKA$^\dag$ and Ning Ning$^\ddag$}
\thanks{\hglue-4.5mm\fontsize{9.6}{9.6}\selectfont$^*$ Department of Mathematics and Computer Science,
	Queensborough Community College, City University of New York
	\\Email: wjliu@qcc.cuny.edu\\$^\dag$ Department of Statistics and Applied Probability, University
	of California, Santa Barbara\\Email: rao@pstat.ucsb.edu\\$^\ddag$ Department of Statistics and Applied Probability, University
	of California, Santa Barbara\\Email: ning@pstat.ucsb.edu\vspace{2mm}}

\maketitle

\vspace{2mm}

\begin{quote}
{\footnotesize \textbf{Abstract:} } \footnotesize It is well known that reconstruction problems, as the interdisciplinary subject, have been studied in numerous
contexts including statistical physics, information theory and computational biology, to name a few. We consider a $2q$-state symmetric model, with two categories of $q$ states in each category, and 3 transition probabilities: the probability to remain in the same state, the probability to change states but remain in the same category, and the probability to change categories.
We construct a nonlinear second order
dynamical system based on this model and show that
the
Kesten-Stigum reconstruction bound is not tight when $q \geq 4$.
\end{quote}

\begin{tabular}{lp{13cm}}
\multicolumn{2}{l} {\footnotesize\it AMS Subject Classification: }
\footnotesize
primary---60K35; 
secondary---82B26; 82B20 
\\
{\footnotesize\it Key words and phrases: } \footnotesize
Kesten-Stigum reconstruction bound; Markov random fields on trees;
Distributional recursion; Dynamical system
\end{tabular}

\section{Introduction}
\label{intro}
\subsection{Preliminaries}
\indent We start with the following broadcasting process that stands as a
discrete, irreducible, aperiodic, and reversible Markov chain. Let
$\mathbb{T}=(\mathbb{V}, \mathbb{E}, \rho)$ be a tree with nodes
$\mathbb{V}$, edges $\mathbb{E}$ and root $\rho\in \mathbb{V}$. Each
edge of the tree acts as a channel on a finite characters set
$\mathcal{C}$, whose elements are configurations on $\mathbb{T}$,
denoted by $\sigma$. We set a probability transition matrix
$\mathbf{M}=(M_{ij})$ as the noisy communication channel on each
edge. The state of the root $\rho$, denoted by $\sigma_\rho$, is
chosen according to an initial distribution $\pi$ on $\mathcal{C}$, and then propagated in the tree as follows: for each
vertex $v$ having $u$ as its parent , the spin at $v$ is defined
according to the probabilities
$$
\mathbf{P}(\sigma_v=j\mid\sigma_u=i)=M_{i j}
$$
with $i, j \in \mathcal{C}$. Roughly speaking, 
reconstruction is to answer the question that considering all the symbols received
at the vertices of the $n$th generation, does this configuration
contain non-vanishing information transmitted by the
root, as $n$ goes to $\infty$?

In this paper, we will restrict our attention to $d$-ary
trees, i.e. the infinite rooted tree where every vertex has
exactly $d$ offspring (every vertex has degree $d+1$ except the root
which has degree $d$). Let $\sigma(n)$ denote the spins at distance
$n$ from the root and let $\sigma^i(n)$ denote $\sigma(n)$
conditioned on $\sigma_\rho = i$. Consider a characters set
$\mathcal{C}=\mathcal{C}_1\cup \mathcal{C}_2$, consisting of two
categories $\mathcal{C}_1=\{1, \ldots, q\}$ and
$\mathcal{C}_2=\{q+1, \ldots, 2q\}$ with $q\geq2$, and the state of
the root $\rho$ is chosen according to the uniform distribution on
$\mathcal{C}$. Moreover, a $2q\times 2q$ probability transition
matrix $\mathbf{M}=(M_{ij})_{2q\times 2q}$ is defined as follows:
\begin{equation*}
M_{ij}= \left\{\begin{array}{ll} p_0 & \quad\textrm{if}\ i=j,
\\
\
\\
p_1& \quad \textrm{if}\ i\neq j\ \textrm{and}\ i, j\ \textrm{are in
the same category},
\\
\
\\
p_2 & \quad \textrm{if}\ i\neq j\ \textrm{and}\ i, j\ \textrm{are in
different categories},
\end{array}
\right.
\end{equation*}
where $p_0$, $p_1$ and $p_2$ are all nonnegative, such that
$p_0+(q-1)p_1+qp_2=1$. It can be verified that the eigenvalues of
$\mathbf{M}$ are $\lambda_1=p_0-p_1$, $\lambda_2=p_0+(q-1)p_1-qp_2$,
and $\lambda_3=p_0+(q-1)p_1+qp_2=1$. Therefore we have two
candidates $\lambda_1$ and $\lambda_2$ for $\lambda$, the second
largest eigenvalue in absolute value, which plays a crucial role in
the reconstruction problem. We now give a formal
definition of the reconstruction.
\begin{definition}
The reconstruction problem for the infinite tree $\mathbb{T}$ is
\emph{solvable} if for some $i, j\in \mathcal{C}$,
$$
\limsup_{n\to \infty}d_{TV}(\sigma^i(n), \sigma^j(n))>0
$$
where $d_{TV}$ is the total variation distance. When the $\limsup$
is $0$, we say the model has \emph{non-reconstruction} on
$\mathbb{T}$.
\end{definition}

\subsection{Background}
\indent Beyond the basic interest in determining the reconstruction
threshold of a Markov random field in probability, this problem is
relevant to statistical physics, biology (Daskalakis et al. \cite{DA}, Mossel \cite{MO2}), and information theory (Bhamidi et al. \cite{BH}, Evans et al. \cite{EV}), 
where one is interested in computing the information capacity of the tree
network.
Most closely related to the origins of this work,
for spin systems
in statistical physics, the threshold for reconstruction is
equivalent to the threshold for extremality of the infinite-volume
Gibbs measure induced by free-boundary conditions, see Georgh \cite{GE}. The
reconstruction threshold also has an important effect
in the efficiency of the Glauber dynamics on trees and random
graphs. It is well known that when the model is reconstructible, the mixing time for the
Glauber dynamics on trees is $n^{1+\Omega(1)}$, while it is slower than at higher temperature when the mixing time is $O(n \log n)$. The corresponding bound is tight for the Ising model, namely,
the mixing time is $O(n \log n)$ when $d\lambda^2< 1$. In Martinelli et al. \cite{MA},
this result is extended to the log Sobolev constant and it is also
shown that for measures on trees, a super-linear decay of
point-to-set correlations implies an $\Omega(1)$ spectral gap for
the Glauber dynamics with free boundary conditions. A similar
transition takes place in the colouring model as shown in Tetali et al. \cite{T}.
Sharp bounds of this type are not known for the hardcore model,
although it is conjectured that the Glauber dynamics should again be
$O(n\log n)$ in the non-reconstruction regime.

For any channel $\mathbf{M}$, it is well known that the
reconstruction problem is connected closely to $\lambda$, the second
largest eigenvalue in absolute value of $\mathbf{M}$. An important
general bound was obtained by Kesten and Stigum ~\cite{KS1,KS2}:
the reconstruction problem is solvable if $d|\lambda|^2>1$
($\lambda$ may be a complex number), which is known as the
Kesten-Stigum bound. On the other hand, for larger noise
($d|\lambda|^2 < 1$) one may wonder whether reconstruction is
possible, by exploiting the whole set of symbols received at the $n$th
generation, through a clever use of the correlations between the
symbols received on the leaves. The answer depends on the channel.

For the binary symmetric channel, it was shown in Bleher et al. \cite{BL} that the
reconstruction problem is solvable if and only if $d\lambda^2>1$.
For all other channels, however, it would be a little challenging to
prove the non-reconstructibility. Mossel~\cite{MO1,MO3} showed that
the Kesten-Stigum bound is not the bound for reconstruction in the
binary asymmetric model with sufficiently large asymmetry or in the
Potts model with sufficiently many characters, which sheds the light on
exploring the tightness of the Kesten-Stigum bound. The first exact
reconstruction threshold in roughly a decade, was obtained by Borgs
et al.~\cite{BO}, in which the authors displayed a delicate analysis of the moment recursion on a weighted
version of the magnetization, and thus achieved a breakthrough result.

A particularly important example is provided by $q$-state symmetric
channels, i.e. Potts models in the terminology of statistical mechanics, with the transition matrix
$$
\mathbf{M}=\left(
  \begin{array}{cccccccc}
    p_0 & p_1 & \cdots & p_1  \\
    p_1 & p_0 & \cdots & p_1  \\
    \vdots & \vdots & \ddots & \vdots \\
    p_1 & p_1 & \cdots & p_0 \\
  \end{array}
\right)
$$
and $\lambda=p_0-p_1$. This model was completely investigated
by Sly~\cite{S} by means of the recursive structure of the tree, and
more importantly, Sly showed 
that non-reconstruction is equivalent to $\lim_{n\to
\infty}x_n=0$, where $x_n=\mathbf{E}\mathbf{P}(\sigma_\rho=1\mid
\sigma^1(n))-\frac1q$. Thus the key idea is to analyze the
recursion relationship between $x_n$ and $x_{n+1}$. This work then goes on to engage the refined recursive equations
of vector-valued distributions and concentration analyses to confirm much of the picture conjectured earlier by M\'{e}zard and Montanari~\cite{MM}.

Inspired by the popular K80 model proposed by Kimura \cite{Kimura1980}, which distinguishes between transitions and transversions, we analyze the case that transition matrix has two mutation classes and $q$ states in each class. Improved flexibility comes along with increased complexity, which is mainly due to the fact that the additional class of mutation complicates the discussion of the second
largest eigenvalue in absolute value. However, by introducing
additional auxiliary quantities $y_n$ and $z_n$ besides $x_n$ defined in Section \ref{sec:Notations}, we succeed in investigating the tightness of the Kesten-Stigum bound.

\subsection{Applications}
\indent The reconstruction problem arises naturally in many fields including statistical physics, where the Ising model and the Potts model are popular  and have been studied extensively from different angles, see \cite{Baxter2003, Cassandro2005, Derrida1994, Dhar2013, Ferrari2002, Giuliani2013, Giuliani2016, Maes1995, Maes1999, Majumdar1992, Olejarz2013, Saul1993, Sire1995a,
Sire1995b, Spirin2001, Spohn1985, Tracy1988, Weeks1973}. In this article, we focus on the reconstruction threshold on trees, which plays an important role in the dynamic phase transitions in certain glassy systems subject to random constraints. For random colorings on the Erd{\"o}s R{\'e}nyi random graph with average connectivity $d$, Achlioptas and Coja-Oghlan \cite{Achlioptas2008} proved that there is a phase transition, from the situation that most of the mass is contained in one giant component, to the case that the space of solutions breaks into exponentially many smaller clusters. This phase transition has been proved corresponding to known bounds on the reconstruction threshold for proper colorings on trees, see e.g. Mossel and Peres \cite{Mossel2003}, Semerjian \cite{Semerjian2008} and Sly \cite{Sly2009}.

In computational biology, the broadcast model is the main model for the evolution
of base pairs of DNA. Phylogenetic reconstruction is a major task of systematic biology, which is to construct the ancestry tree of a collection of species, given the information of present species. The corresponding reconstruction threshold answers the question whether the ancestral DNA information can be reconstructed from a known phylogenetic tree. This threshold is also crucial to determine the number of samples required, in the sense that, only enumerations of each type of spin at the leaves are collected, regardless of their positions on the leaves. Interested readers on Phylogenetic tree reconstruction are referred to Roch \cite{Roch2006} and Daskalakis et al. \cite{DA}.

The popular K80 model \cite{Kimura1980}, has some obvious advantages over other models in Phylogeny reconstruction, which is favored by both Akaike Information Criterion and Bayesian Information Criterion (see Section  2.2.2 in Cadotte and Davies \cite{Cadotte2016}). The K80 model distinguishes between transitions ($A \leftrightarrow G$, i.e. from purine to purine, or $C \leftrightarrow T$, i.e. from pyrimidine to pyrimidine) and transversions (from purine to pyrimidine or vice versa). Inspired by this and related literatures, we analyze the case that the transition matrix has two mutation classes and $q$ states in each class. We believe that the $q$-state symmetric Potts model as a generalization of $2$-state symmetric Ising model, cannot fully represent the spirit of the classical $2$-state symmetric Ising model in terms of dichotomy, and this is one of the areas this work can contribute to.

A tree is a connected undirected graph with no simple circuits.
In other words, an undirected graph is a tree if and only if there is a unique simple path between any two of its
vertices. The theory that the reconstruction threshold on trees corresponds to the reconstruction threshold on locally treelike graphs, is verified in Gerschenfeld and Montanari \cite{Gerschenfeld2007}. The strong and increasing interest in the study of the properties of social networks, is a result of the rapid and global emergence of online social networks and their meteoric adoption by millions of Internet users. When it comes to Socio--psychological mechanisms of generation and dissemination of network, our model's advantage in providing more flexibility to mimic psychological behaviors is obvious. For example, our model and the construction threshold can be used to effectively identify community effect in social networks and customer loyalty in marketing research, especially for different firms or organizations who want to promote their products or philosophies.
In this sense, many possible extensions can be made on research on graph structures with psychological factors involved, such as the work by Liu, Ying and Shakkottai \cite{Liu2010} on Ising model based analysis on the formation and propagation of opinions across networks,  the work by Bisconti et al. \cite{Bisconti2015} on Potts model based analysis on the reconstruction of a real world social network and loopy belief propagation, etc.

\subsection{Main Results and Proof Sketch}
\indent Because non-reconstruction happens at most $d|\lambda|^2=1$, without
loss of generality, it would be convenient to presume $1/2\leq
d|\lambda|^2\leq 1$ in the following context.
\begin{Main Theorem}
\label{qgeq4} Assume $0<|\lambda_2|\leq |\lambda_1|$. When $q\geq
4$, for every $d$ the Kesten-Stigum bound is not tight, i.e. the
reconstruction is solvable for some $\lambda_1$ even if
$d\lambda_1^2<1$.
\end{Main Theorem}

The ideas and techniques used to prove the Main Theorem can be seen as the following.
One standard to classify reconstruction and nonreconstruction is to analyze the quantity $x_n$: the probability of giving a correct guess of the root given the spins $\sigma(n)$ at distance $n$ from the root, minus the probability of guessing the root randomly which is $\frac{1}{2q}$ in this case. Nonreconstruction means that the mutual information between the root and the spins at distance $n$ goes to $0$ as $n$ tends to infinity. It can be established that the nonreconstruction is equivalent to
$$
\lim_{n\to \infty}x_n=0.
$$

Our analysis is similar to Borgs et al. \cite{BO}, Chayes et al. \cite{Chayes1986} in the context of spin--glasses, and Sly \cite{S}. However, the two classes of mutation complicates the discussion of $\lambda$, the second largest eigenvalue in absolute value of the transition matrix, which makes the problem much more challenging. In this case, it is necessary to consider the corresponding quantities similar to $x_n$, viz. wrong guess but right group $y_n$, and wrong guess and even wrong group $z_n$. In Section 2.2, we investigate the properties and relation between $x_n$,
$y_n$ and $z_n$. By these preliminary results, we
focus on the analysis of $x_n$ and $z_n$ in the sequel.

\indent In order to research the reconstruction, according to the Markov random field property, we establish the distributional recursion and moment recursion, by analyzing the recursive relation between the $n$th and the $(n+1)$th generations' structure of the tree. Furthermore, we display that the interactions between spins become very weak, if they are sufficiently far away from each other. Therefore, we can obtain a nonlinear dynamical system. If $x_n$ is small, we are able to develop the concentration analysis and achieve the approximation to the dynamical system:
\begin{equation*}
\label{X}
\left\{\begin{array}{ll}
x_{n+1}\approx d\lambda_1^2x_n+(d\lambda_1^2-d\lambda_2^2)z_n+\frac{d(d-1)}{2}\left(\frac{q(2q-5)}{q-1}\lambda_1^4(x_n+z_n)^2-4q\lambda_1^2\lambda_2^2(x_n+z_n)z_n-4q\lambda_2^4z_n^2\right)
\\
\
\\
z_{n+1}\approx d\lambda_2^2z_n-\frac{d(d-1)}{2}\left(\frac{q}{q-1}\lambda_1^4(x_n+z_n)^2-4q\lambda_2^4z_n^2\right).
\end{array}
\right.
\end{equation*}

 Finally, we investigate the stability of the system. We establish the threshold of $q$ relevant to the reconstruction. When $q \geq 4$, even if $d\lambda_1^2<1$ for some $\lambda_1$, $x_n$ will not converge to $0$ and hence there is reconstruction beyond the Kesten-Stigum bound. More detailed definitions and interpretations can be seen in the next Section.


\section{SECOND ORDER RECURSION RELATION}
\subsection{Notations}
\label{sec:Notations}
\indent Let $u_1,\ldots,u_d$ be the children of $\rho$ and $\mathbb{T}_v$ be
the subtree of descendants of $v\in \mathbb{T}$. Furthermore, if we
set $d(\cdot, \cdot)$ as the graph-metric distance on $\mathbb{T}$,
denote the $n$th level of the tree by $L_n=\{v\in \mathbb{V}:
d(\rho, v)=n\}$ and then let $\sigma_j(n)$ be the spins on
$L_n\cap \mathbb{T}_{u_j}$. For a configuration $A$ on $L_n$, define
the posterior function
\begin{equation}
\label{posterior}
f_n(i, A)=\textbf{P}(\sigma_\rho=i\mid\sigma(n)=A).
\end{equation}
By the recursive nature of the tree for a configuration $A$ on $L(n
+ 1) \cap \mathbb{T}_{u_j}$, there is an equivalent form
$$
f_n(i, A)=\textbf{P}(\sigma_{u_j}=i\mid\sigma_j(n+1)=A).
$$
Now for any $1\leq i\leq 2q$, define a collection of random
variables
$$
X_i(n)=f_n(i, \sigma(n))
$$
to describe the posterior probability of state $i$ at the root given
the random configuration $\sigma(n)$ of the leaves, and analogously,
$$
X^{(1)}(n)=f_n(1, \sigma^1(n)), \quad X^{(2)}(n)=f_n(2,
\sigma^1(n)),\quad X^{(3)}(n)=f_n(q+1, \sigma^1(n)).
$$
By symmetry, the collections $\{f_n(i, \sigma^1(n)): 2\leq i\leq q\}$
and $\{f_n(i, \sigma^1(n)): q+1\leq i\leq 2q\}$ are exchangeable
respectively; in addition, $f_n(j, \sigma^i(n))$ is distributed as
$$
f_n(j, \sigma^i(n))\stackrel{\mathbb{D}}{\sim}
\left\{\begin{array}{ll} X^{(1)}(n) & \quad \textup{if}\ i=j,
\\
\
\\
X^{(2)}(n)& \quad \textup{if}\ i\neq j\ \textup{are in the same
category},
\\
\
\\
X^{(3)}(n) & \quad \textup{if}\ i\neq j\ \textup{are in different
categories}.
\end{array}
\right.
$$
Finally, denote the first and second central moments of $X^{(1)}(n)$, $X^{(2)}(n)$ and $X^{(3)}(n)$,
which would be the principal quantities in our analysis, as
$$
x_n=\mathbf{E}\left(X^{(1)}(n)-\frac{1}{2q}\right), \quad
y_n=\mathbf{E}\left(X^{(2)}(n)-\frac{1}{2q}\right), \quad
z_n=\mathbf{E}\left(X^{(3)}(n)-\frac{1}{2q}\right),
$$
and
$$
u_n=\mathbf{E}\left(X^{(1)}(n)-\frac{1}{2q}\right)^2,\quad
v_n=\mathbf{E}\left(X^{(2)}(n)-\frac{1}{2q}\right)^2,\quad
w_n=\mathbf{E}\left(X^{(3)}(n)-\frac{1}{2q}\right)^2.
$$

\subsection{Preliminaries}
\label{sec:Preparations}
\indent For any $i=1,\cdots,2q$ and nonnegative $n\in \mathbb{Z}$, it
is concluded from the symmetric property of the tree that
$$\mathbf{E}X_i(n)=\frac1{2q}$$ is always true.

\begin{lemma}
\label{thm:xn} For any $n\in \mathbb{N}\cup\{0\}$, we have
$$
x_n=\mathbf{E}\sum_{i=1}^{2q}\left(X_i(n)-\frac{1}{2q}\right)^2\geq
0, \quad z_n\leq 0, \quad \textup{and}\quad x_n+z_n\geq 0.
$$
\end{lemma}
\begin{proofsect}{Proof}
First, by Bayes' rule, we have
$$
x_n+\frac{1}{2q}=\sum_Af_n(1,
A)\mathbf{P}(\sigma(n)=A\mid\sigma_\rho=1)
=2q\sum_A\mathbf{P}(\sigma(n)=A)f_n^2(1, A) =2q\mathbf{E}X_1^2(n)
$$
and 
\begin{equation}
\label{xn}
0\leq\mathbf{E}\sum_{i=1}^{2q}\left(X_i(n)-\frac{1}{2q}\right)^2=\sum_{i=1}^{2q}\mathbf{E}X_i^2(n)-\frac{2}{2q}\sum_{i=1}^{2q}\mathbf{E}X_i(n)+\frac{1}{2q}=x_n.
\end{equation}

Next, we consider the covariance matrix of random variables
$\left\{X_i(n)-\frac{1}{2q}\right\}_1^{2q}$ and express covariances in terms of $x_n$, $y_n$ and
$z_n$. Similarly, we obtain
$$
y_n+\frac{1}{2q}=2q\sum_A\mathbf{P}(\sigma(n)=A)f_n(1, A)f_n(2, A)
=2q\mathbf{E}X_1(n)X_2(n),
$$
so for any $i_1<i_2$ in the same category, it is concluded from the
symmetric property of the tree that
$$
\mathbf{E}\left(X_{i_1}(n)-\frac{1}{2q}\right)\left(X_{i_2}(n)-\frac{1}{2q}\right)=\mathbf{E}\left(X_1-\frac{1}{2q}\right)\left(X_2-\frac{1}{2q}\right)=\frac{y_n}{2q}.
$$
Similarly, if $i_1$ and $i_2$ are from different categories, we have
$$
\mathbf{E}\left(X_{i_1}(n)-\frac{1}{2q}\right)\left(X_{i_2}(n)-\frac{1}{2q}\right)=\frac{z_n}{2q}.
$$
Therefore, the covariance matrix is given by
$$
\Sigma_X(n)=\left(
  \begin{array}{cccccccc}
    \frac{x_n}{2q} & \frac{y_n}{2q} & \cdots & \frac{y_n}{2q} & \frac{z_n}{2q} & \frac{z_n}{2q} & \cdots & \frac{z_n}{2q} \\
    \frac{y_n}{2q} & \frac{x_n}{2q} & \cdots & \frac{y_n}{2q} & \frac{z_n}{2q} & \frac{z_n}{2q} & \cdots & \frac{z_n}{2q} \\
    \vdots & \vdots & \ddots & \vdots & \vdots & \vdots & \ddots & \vdots \\
    \frac{y_n}{2q} & \frac{y_n}{2q} & \cdots & \frac{x_n}{2q} & \frac{z_n}{2q} & \frac{z_n}{2q} & \cdots & \frac{z_n}{2q} \\
    \frac{z_n}{2q} & \frac{z_n}{2q} & \cdots & \frac{z_n}{2q} & \frac{x_n}{2q} & \frac{y_n}{2q} & \cdots & \frac{y_n}{2q} \\
    \frac{z_n}{2q} & \frac{z_n}{2q} & \cdots & \frac{z_n}{2q} & \frac{y_n}{2q} & \frac{x_n}{2q} & \cdots & \frac{y_n}{2q} \\
    \vdots & \vdots & \ddots & \vdots & \vdots & \vdots & \ddots & \vdots \\
    \frac{z_n}{2q} & \frac{z_n}{2q} & \cdots & \frac{z_n}{2q} & \frac{y_n}{2q} & \frac{y_n}{2q} & \cdots & \frac{x_n}{2q} \\
  \end{array}
\right)_{2q\times 2q}
$$
whose eigenvalues are $0$, $\frac{x_n+(q-1)y_n-qz_n}{2q}$ and
$\frac{x_n-y_n}{2q}$. It is well known that the
covariance matrix of a multivariate probability distribution is
always positive semi-definite, which implies that all eigenvalues
are nonnegative, say, $x_n+(q-1)y_n-qz_n\geq 0$ and $x_n-y_n\geq 0$. It suffices to
complete the proof, by these results and the fact $x_n+(q-1)y_n+qz_n=0$.
\end{proofsect}

\begin{lemma}
\label{xnun} For any $n\in \mathbb{N}\cup\{0\}$, the following hold:
\begin{enumerate}[(i)]
\item $x_n=u_n+(q-1)v_n+qw_n$;

\item $\mathbf{E}\left(X^{(1)}(n)-\frac{1}{2q}\right)\left(X^{(2)}(n)-\frac{1}{2q}\right)=v_n+\frac{y_n-x_n}{2q}$;

\item $\mathbf{E}\left(X^{(1)}(n)-\frac{1}{2q}\right)\left(X^{(3)}(n)-\frac{1}{2q}\right)=w_n+\frac{z_n-x_n}{2q}$;

\item $\mathbf{E}\left(X^{(2)}(n)-\frac{1}{2q}\right)\left(X^{(3)}(n)-\frac{1}{2q}\right)=-\frac{w_n}{q-1}-\frac{z_n}{2q(q-1)}-\frac{y_n}{2q}$;

\item $\mathbf{E}\left(f_n(q+1, \sigma^1(n))-\frac{1}{2q}\right)\left(f_n(2q, \sigma^1(n))-\frac{1}{2q}\right)=-\frac{w_n}{q-1}-\frac{z_n}{2(q-1)}$;

\item $\mathbf{E}\left(f_n(2, \sigma^1(n))-\frac{1}{2q}\right)\left(f_n(q, \sigma^1(n))-\frac{1}{2q}\right)=-\frac{2v_n}{q-2}-\frac{z_n}{2(q-1)}+\frac{qw_n}{(q-1)(q-2)}$.
\end{enumerate}
\end{lemma}
\begin{proofsect}{Proof}
By the total probability formula and using Lemma \ref{thm:xn}, we can prove (i) as follows:
\begin{eqnarray*}
x_n&=&\mathbf{E}\sum_{i=1}^{2q}\left(X_i(n)-\frac{1}{2q}\right)^2
\\
&=&\sum_{j=1}^{2q}\mathbf{E}\left(\sum_{i=1}^{2q}\left(X_i(n)-\frac{1}{2q}\right)^2\mid\sigma_\rho=j\right)\mathbf{P}(\sigma_\rho=j)
\\
&=&\sum_{j=1}^{2q}\frac{1}{2q}\left[\mathbf{E}\left(X^{(1)}(n)-\frac{1}{2q}\right)^2+(q-1)\mathbf{E}\left(X^{(2)}(n)-\frac{1}{2q}\right)^2+q\mathbf{E}\left(X^{(3)}(n)-\frac{1}{2q}\right)^2\right]
\\
&=&\mathbf{E}\left(X^{(1)}(n)-\frac{1}{2q}\right)^2+(q-1)\mathbf{E}\left(X^{(2)}(n)-\frac{1}{2q}\right)^2+q\mathbf{E}\left(X^{(3)}(n)-\frac{1}{2q}\right)^2
\\
&=&u_n+(q-1)v_n+qw_n.
\end{eqnarray*}
	Applying the same technique, we obtain
	\begin{equation*}
	\begin{aligned}
	\mathbf{E}X^{(1)}(n)X^{(2)}(n)&=\sum_{A}\mathbf{P}(\sigma_\rho=1\mid\sigma(n)=A)\mathbf{P}(\sigma_\rho=2\mid\sigma(n)=A)P(\sigma(n)=A\mid\sigma_\rho=1)
	\\
	&=\sum_{A}[\mathbf{P}(\sigma_\rho=1\mid\sigma(n)=A)]^2\mathbf{P}(\sigma(n)=A\mid\sigma_\rho=2)
	\\
	&=\mathbf{E}\left(X^{(2)}(n)\right)^2
	\end{aligned}
	\end{equation*}
	and hence (ii) follows:
	\begin{eqnarray*}
		\mathbf{E}\left(X^{(1)}(n)-\frac{1}{2q}\right)\left(X^{(2)}(n)-\frac{1}{2q}\right)&=&\mathbf{E}\left(X^{(2)}-\frac{1}{2q}\right)^2+\frac{y_n-x_n}{2q}
		=v_n+\frac{y_n-x_n}{2q}.
	\end{eqnarray*}
	Similarly, (iii) turns out to be true due to
	$$\mathbf{E}X^{(1)}(n)X^{(3)}(n)=\mathbf{E}\left(X^{(3)}(n)\right)^2.$$
	The statement
	(iv), (v) and (vi) can be handled in the same way, using the
	symmetry,
	\begin{equation}
	\label{X2X3}
	\begin{aligned}
	\mathbf{E}X^{(2)}(n)X^{(3)}(n)
	&=\sum_{A}\mathbf{P}(\sigma_\rho=2\mid\sigma(n)=A)\mathbf{P}(\sigma_\rho=q+1\mid\sigma(n)=A)\mathbf{P}(\sigma(n)=A\mid\sigma_\rho=1)
	\\
	&=\sum_{A}\mathbf{P}(\sigma_\rho=1\mid\sigma(n)=A)\mathbf{P}(\sigma_\rho=2\mid\sigma(n)=A)\mathbf{P}(\sigma(n)=A\mid\sigma_\rho=q+1)
	\\
	&=\sum_{A}\mathbf{P}(\sigma_\rho=q+1\mid\sigma(n)=A)\mathbf{P}(\sigma_\rho=2q\mid\sigma(n)=A)\mathbf{P}(\sigma(n)=A\mid\sigma_\rho=1)
	\\
	&=\mathbf{E}f_n(q+1, \sigma^1(n))f_n(2q, \sigma^1(n)).
	\end{aligned}
	\end{equation}
	To obtain $\mathbf{E}X^{(2)}(n)X^{(3)}(n)$, recall that
	\begin{eqnarray*}
		z_n+\frac{1}{2q}&=&\mathbf{E}X^{(3)}(n)
		\\
		&=&\mathbf{E}f_n(q+1, \sigma^1(n))\sum_{i=1}^{2q}f_n(i, \sigma^1(n))
		\\
		&=&\mathbf{E}X^{(1)}(n)X^{(3)}(n)+(q-1)\mathbf{E}X^{(2)}(n)X^{(3)}(n)+\mathbf{E}(X^{(3)})^2
		+(q-1)\mathbf{E}f_n(q+1, \sigma^1(n))f_n(2q, \sigma^1(n))
		\\
		&=&\mathbf{E}X^{(1)}(n)X^{(3)}(n)+2(q-1)\mathbf{E}X^{(2)}(n)X^{(3)}(n)+\mathbf{E}(X^{(3)})^2,
	\end{eqnarray*}
	which implies that
	\begin{equation}
	\label{X2X3expression}
	\mathbf{E}\left(X^{(2)}(n)-\frac{1}{2q}\right)\left(X^{(3)}(n)-\frac{1}{2q}\right)
	=-\frac{w_n}{q-1}-\frac{z_n}{2q(q-1)}-\frac{y_n}{2q}.
	\end{equation}
	Thus,~\eqref{X2X3} together with~\eqref{X2X3expression} gives
	\begin{equation*}
	\mathbf{E}\left(f_n(q+1,
	\sigma^1(n))-\frac{1}{2q}\right)\left(f_n(2q,
	\sigma^1(n))-\frac{1}{2q}\right)
	=-\frac{w_n}{q-1}-\frac{z_n}{2(q-1)}.
	\end{equation*}
	As in the preceding discussion, considering
	\begin{eqnarray*}
		y_n+\frac{1}{2q}&=&\mathbf{E}f_n(2, \sigma^1(n))\sum_{i=1}^{2q}f(i,
		\sigma^1(n))
		\\
		&=&2\mathbf{E}(X^{(2)}(n))^2+(q-2)\mathbf{E}f_n(2, \sigma^1(n))f_n(q,
		\sigma^1(n))+q\mathbf{E}X^{(2)}(n)X^{(3)}(n),
	\end{eqnarray*}
	we have
	\begin{equation*}
	\begin{aligned}
	&\mathbf{E}\left(f_n(2, \sigma^1(n))-\frac{1}{2q}\right)\left(f_n(q,
	\sigma^1(n))-\frac{1}{2q}\right)
	\\
	&=\frac{1}{q-2}\left(y_n+\frac{1}{2q}-2\mathbf{E}(X^{(2)}(n))^2-q\mathbf{E}X^{(2)}(n)X^{(3)}(n)\right)-\frac{2}{2q}\left(y_n+\frac{1}{2q}\right)+\frac{1}{4q^2}
	\\
	&=-\frac{2v_n}{q-2}-\frac{z_n}{2(q-1)}+\frac{qw_n}{(q-1)(q-2)}.
	\end{aligned}
	\end{equation*}
\end{proofsect}

\subsection{Means and Covariances of $Y_{ij}$}
\label{Yij} 
\indent Defining
$$
Y_{ij}(n) = f_n\left(i, \sigma^1_j (n + 1)\right),
$$
and taking advantage of the symmetries of the model, it is apparent
that the random vectors $(Y_{ij})_{i=1}^{2q}$
are independent, for $j=1, \ldots, d$. The central moments of $Y_{ij}$ would play a key role
in further analysis, therefore it is necessary to figure them out in the
first place. For each $1\leq j\leq d$, we rely on the total
probability formula to conclude:
\begin{enumerate}[(i)]
\item when $i=1$,
\begin{equation*}
\begin{aligned}
\mathbf{E}\left(Y_{1j}(n)-\frac{1}{2q}\right)
&=p_0\mathbf{E}\left(X^{(1)}(n)-\frac{1}{2q}\right)+(q-1)p_1\mathbf{E}\left(X^{(2)}(n)-\frac{1}{2q}\right)+qp_2\mathbf{E}\left(X^{(3)}(n)-\frac{1}{2q}\right)
\\
&=\lambda_1x_n+(\lambda_1-\lambda_2)z_n;
\end{aligned}
\end{equation*}

\item for $2\leq i\leq q$,
\begin{equation*}
\begin{aligned}
\mathbf{E}\left(Y_{ij}(n)-\frac{1}{2q}\right)
&=p_1\mathbf{E}\left(X^{(1)}(n)-\frac{1}{2q}\right)+[p_0+(q-2)p_1]\mathbf{E}\left(X^{(2)}(n)-\frac{1}{2q}\right)+qp_2\mathbf{E}\left(X^{(3)}(n)-\frac{1}{2q}\right)
\\
&=-\frac{\lambda_1}{q-1}x_n-\frac{\lambda_1+(q-1)\lambda_2}{q-1}z_n;
\end{aligned}
\end{equation*}

\item for $q+1\leq i\leq 2q$, by means of the identity $\sum_{i=1}^{2q}Y_{ij}(n)\equiv1$, it follows
immediately that
$$
\mathbf{E}\left(Y_{ij}(n)-\frac{1}{2q}\right)=-\frac{1}{q}\sum_{i=1}^q\mathbf{E}\left(Y_{ij}(n)-\frac{1}{2q}\right)=\lambda_2z_n;
$$

\item resembling the discussion of (i), (ii) and (iii), it is further concluded that when $i=1$,
\begin{equation*}
\mathbf{E}\left(Y_{1j}(n)-\frac{1}{2q}\right)^2
=\frac{1+\lambda_2-2\lambda_1}{2q}x_n+\lambda_1u_n+(\lambda_1-\lambda_2)w_n;
\end{equation*}

\item for $2\leq i\leq q$,
\begin{equation*}
\mathbf{E}\left(Y_{ij}(n)-\frac{1}{2q}\right)^2
=\left(\frac{1}{2q}+\frac{\lambda_2}{2q}+\frac{\lambda_1}{q(q-1)}\right)x_n-\frac{\lambda_1}{q-1}u_n-\frac{\lambda_1+(q-1)\lambda_2}{q-1}w_n;
\end{equation*}

\item for $q+1\leq i\leq 2q$,
\begin{equation*}
\mathbf{E}\left(Y_{ij}(n)-\frac{1}{2q}\right)^2
=\frac{1-\lambda_2}{2q}x_n+\lambda_2w_n;
\end{equation*}

\item for $2\leq i\leq q$,
\begin{equation*}
\begin{aligned}
\mathbf{E}\left(Y_{1j}(n)-\frac{1}{2q}\right)\left(Y_{ij}(n)-\frac{1}{2q}\right)
=\frac{(q+2)\lambda_1-\lambda_2-1}{2q(q-1)}x_n-\frac{z_n}{2(q-1)}-\frac{\lambda_1}{q-1}u_n-\frac{(q+1)\lambda_1-\lambda_2}{q-1}w_n;
\end{aligned}
\end{equation*}

\item for $q+1\leq i\leq 2q$,
\begin{equation*}
\mathbf{E}\left(Y_{1j}(n)-\frac{1}{2q}\right)\left(Y_{ij}(n)-\frac{1}{2q}\right)
=-\frac{\lambda_1}{2q}x_n+\frac{z_n}{2q}+\lambda_1w_n;
\end{equation*}

\item for $1< i_1< i_2\leq q$,
\begin{equation*}
\label{recursion}
\begin{aligned}
&\mathbf{E}\left(Y_{i_1j}(n)-\frac{1}{2q}\right)\left(Y_{i_2j}(n)-\frac{1}{2q}\right)
\\
&=\left[-\frac{2(q+2)\lambda_1+(q-2)\lambda_2}{2q(q-1)(q-2)}-\frac{1}{2q(q-1)}\right]x_n-\frac{z_n}{2(q-1)}
+\frac{2\lambda_1}{(q-1)(q-2)}u_n+\frac{2(q+1)\lambda_1+(q-2)\lambda_2}{(q-1)(q-2)}w_n;
\end{aligned}
\end{equation*}

\item for $1< i_1\leq q< i_2\leq 2q$,
\begin{eqnarray*}
\mathbf{E}\left(Y_{i_1j}(n)-\frac{1}{2q}\right)\left(Y_{i_2j}(n)-\frac{1}{2q}\right)
&=&\frac{\lambda_1}{2q(q-1)}x_n+\frac{z_n}{2q}-\frac{\lambda_1}{q-1}w_n;
\end{eqnarray*}

\item for $q+1\leq i_1< i_2\leq 2q$,
\begin{eqnarray*}
\mathbf{E}\left(Y_{i_1j}(n)-\frac{1}{2q}\right)\left(Y_{i_2j}(n)-\frac{1}{2q}\right)
=\frac{\lambda_2-1}{2q(q-1)}x_n-\frac{z_n}{2(q-1)}-\frac{\lambda_2}{q-1}w_n.
\end{eqnarray*}
\end{enumerate}

\subsection{Distributional Recursion}
\label{Zresults} 
\indent The key method of this paper is to analyze the relation between
 $X^{(1)}(n), X^{(3)}(n)$ and $X^{(1)}(n+1), X^{(3)}(n+1)$ using the recursive structure of the tree. Take $A=\sigma^1(n+1)$ and then the following
relation follows from the Markov random field property:
$$
X^{(1)}(n+1)=f_{n+1}(1,
\sigma^1(n+1))=\frac{Z_1}{\sum_{i=1}^{2q}Z_i}
$$
and
$$
X^{(3)}(n+1)=f_{n+1}(q+1,
\sigma^1(n+1))=\frac{Z_{q+1}}{\sum_{i=1}^{2q}Z_i},
$$
where
\begin{enumerate}[(A)]
\item for $1\leq i\leq q$,
\begin{eqnarray*}
Z_i=Z_i(n)&=&\prod_{j=1}^d\left[1+2q(p_0-p_2)\left(Y_{ij}-\frac{1}{2q}\right)+2q(p_1-p_2)\sum_{1\leq\ell\neq
i\leq q}\left(Y_{\ell j}-\frac{1}{2q}\right)\right]
\\
&=&\prod_{j=1}^d\left[1+2q(p_0-p_1)\left(Y_{ij}-\frac{1}{2q}\right)-2q(p_1-p_2)\sum_{q+1\leq\ell\leq
2q}\left(Y_{\ell j}-\frac{1}{2q}\right)\right]
\\
&=&\prod_{j=1}^d\left[1+2q\lambda_1\left(Y_{ij}-\frac{1}{2q}\right)+2(\lambda_1-\lambda_2)\sum_{q+1\leq\ell\leq
2q}\left(Y_{\ell j}-\frac{1}{2q}\right)\right]
\end{eqnarray*}

\item for $q+1\leq i\leq 2q$,
\begin{eqnarray*}
Z_i=Z_i(n)&=&\prod_{j=1}^d\left[1+2q(p_0-p_2)\left(Y_{ij}-\frac{1}{2q}\right)+2q(p_1-p_2)\sum_{q+1\leq\ell\neq
i\leq 2q}\left(Y_{\ell j}-\frac{1}{2q}\right)\right]
\\
&=&\prod_{j=1}^d\left[1+2q(p_0-p_1)\left(Y_{ij}-\frac{1}{2q}\right)-2q(p_1-p_2)\sum_{1\leq\ell\leq
q}\left(Y_{\ell j}-\frac{1}{2q}\right)\right]
\\
&=&\prod_{j=1}^d\left[1+2q\lambda_1\left(Y_{ij}-\frac{1}{2q}\right)+2(\lambda_1-\lambda_2)\sum_{1\leq\ell\leq
q}\left(Y_{\ell j}-\frac{1}{2q}\right)\right].
\end{eqnarray*}
\end{enumerate}

To continue the proof, it is necessary to firstly derive some
identities concerning $Z_i(n)$.
\begin{lemma}
\label{Z1Zi} For any nonnegative $n\in \mathbb{Z}$ and $1\leq i\leq
2q$, we have
$$
\mathbf{E}Z_1(n) Z_i(n)=\mathbf{E}Z_i(n)^2,
$$
and given any $2\leq i_1\leq q<q+1\leq i_2\leq 2q$, we have
$$
\mathbf{E}Z_{i_1}(n)Z_{i_2}(n)=\mathbf{E}Z_{q+1}(n)Z_{2q}(n).
$$
\end{lemma}
\begin{proofsect}{Proof}
When $i=1$, the result is trivial. If $2\leq i\leq 2q$, for any
configurations $A=(A_1, \ldots, A_d)$ on the $(n+1)$th level, where
$A_j$ denote the spins on $L_{n+1}\cap \mathbb{T}_{u_j}$ ,we have
\begin{eqnarray*}
Z_i(A)=2q\frac{\mathbf{P}(\sigma(n+1)=A)}{\prod_{j=1}^d\mathbf{P}(\sigma_j(n+1)=A_j)}\mathbf{P}(\sigma_\rho=i\mid\sigma(n+1)=A)
\end{eqnarray*}
By the symmetry of the tree, we have
\begin{equation*}
\begin{aligned}
\mathbf{E}Z_1Z_i&=(2q)^2\sum_A\left(\frac{\mathbf{P}(\sigma(n+1)=A)}{\prod_{j=1}^d\mathbf{P}(\sigma_j(n+1)=A_j)}\right)^2\mathbf{P}(\sigma_\rho=1\mid\sigma(n+1)=A)
\\
&\quad\times\mathbf{P}(\sigma_\rho=i\mid\sigma(n+1)=A)\mathbf{P}(\sigma(n+1)=A\mid\sigma_\rho=1)
\\
&=(2q)^2\sum_A\left(\frac{\mathbf{P}(\sigma(n+1)=A)}{\prod_{j=1}^d\mathbf{P}(\sigma_j(n+1)=A_j)}\right)^2\mathbf{P}^2(\sigma_\rho=1\mid\sigma(n+1)=A)
\times\mathbf{P}(\sigma(n+1)=A\mid\sigma_\rho=i)
\\
&=(2q)^2\sum_A\left(\frac{\mathbf{P}(\sigma(n+1)=A)}{\prod_{j=1}^d\mathbf{P}(\sigma_j(n+1)=A_j)}\right)^2\mathbf{P}^2(\sigma_\rho=i\mid\sigma(n+1)=A)
\mathbf{P}(\sigma(n+1)=A\mid\sigma_\rho=1)
\\
&=\mathbf{E}Z_i^2.
\end{aligned}
\end{equation*}
Similarly, given arbitrary $2\leq i_1\leq q< i_2\leq 2q$, there is
\begin{equation*}
\begin{aligned}
\mathbf{E}Z_{i_1}Z_{i_2}&=(2q)^2\sum_A\left(\frac{\mathbf{P}(\sigma(n+1)=A)}{\prod_{j=1}^d\mathbf{P}(\sigma_j(n+1)=A_j)}\right)^2\mathbf{P}(\sigma_\rho=1\mid\sigma(n+1)=A)
\\
&\quad\times\mathbf{P}(\sigma_\rho=i_1\mid\sigma(n+1)=A)\mathbf{P}(\sigma(n+1)=A\mid\sigma_\rho=i_2)
\\
&=(2q)^2\sum_A\left(\frac{\mathbf{P}(\sigma(n+1)=A)}{\prod_{j=1}^d\mathbf{P}(\sigma_j(n+1)=A_j)}\right)^2\mathbf{P}(\sigma_\rho=q+1\mid\sigma(n+1)=A)
\\
&\quad\times\mathbf{P}(\sigma_\rho=2q\mid\sigma(n+1)=A)\mathbf{P}(\sigma(n+1)=A\mid\sigma_\rho=1)
\\
&=\mathbf{E}Z_{q+1}Z_{2q}.
\end{aligned}
\end{equation*}
\end{proofsect}

We now calculate approximations for the means and
variances of monomials of the $Z_i$ by expanding them using Taylor
series, similar to Lemma
2.6~\cite{S}. In the following, note that the $O_q$ terms depend only on $q$.
\begin{enumerate}[(i)]
\item When $i=1$,
\begin{eqnarray*}
\mathbf{E}Z_1
=1+d\left[2q\lambda_1^2x_n+2q(\lambda_1^2-\lambda_2^2)z_n\right]+\frac{d(d-1)}{2}\left[2q\lambda_1^2x_n+2q(\lambda_1^2-\lambda_2^2)z_n\right]^2+O_q(x_n^3);
\end{eqnarray*}

\item For $2\leq i\leq q$,
\begin{eqnarray*}
\mathbf{E}Z_i=1+d\left[-\frac{2q\lambda_1^2}{q-1}x_n-\left(\frac{2q\lambda_1^2}{q-1}+2q\lambda_2^2\right)z_n\right]+\frac{d(d-1)}{2}\left[-\frac{2q\lambda_1^2}{q-1}x_n-\left(\frac{2q\lambda_1^2}{q-1}+2q\lambda_2^2\right)z_n\right]^2+O_q(x_n^3);
\end{eqnarray*}

\item For $q+1\leq i\leq 2q$,
\begin{eqnarray*}
\mathbf{E}Z_i=1+d(2q\lambda_2^2z_n)+\frac{d(d-1)}{2}(2q\lambda_2^2z_n)^2+O_q(x_n^3).
\end{eqnarray*}
\end{enumerate}

Next consider covariances of $Z_i$s. By Lemma~\ref{Z1Zi} it is
known that $\mathbf{E}Z_1Z_i=\mathbf{E}Z_i^2$, so that we can skip
calculating these terms:
\begin{enumerate}[(a)]
\item when $i=1$,
\begin{eqnarray*}
\mathbf{E}Z_1^2=1+d\Pi_1+\frac{d(d-1)}{2}\Pi_1^2+O_q(x_n^3),
\end{eqnarray*}
where
\begin{eqnarray*}
\Pi_1&=&\mathbf{E}\left[1+2q\lambda_1\left(Y_{1j}-\frac{1}{2q}\right)+2(\lambda_1-\lambda_2)\sum_{q+1\leq\ell\leq
2q}\left(Y_{\ell j}-\frac{1}{2q}\right)\right]^2-1
\\
&=&6q\lambda_1^2x_n+6q(\lambda_1^2-\lambda_2^2)z_n+4q^2\lambda_1^3\left(u_n-\frac{x_n}{2q}\right)
+12q^2\lambda_1^2(\lambda_1-\lambda_2)\left(w_n-\frac{x_n}{2q}\right);
\end{eqnarray*}

\item for $2\leq i\leq q$,
\begin{eqnarray*}
\mathbf{E}Z_i^2=1+d\Pi_2+\frac{d(d-1)}{2}\Pi_2^2+O_q(x_n^3)
\end{eqnarray*}
where
\begin{eqnarray*}
\Pi_2&=&\mathbf{E}\left[1+2q\lambda_1\left(Y_{ij}-\frac{1}{2q}\right)+2(\lambda_1-\lambda_2)\sum_{q+1\leq\ell\leq
2q}\left(Y_{\ell j}-\frac{1}{2q}\right)\right]^2-1
\\
&=&\frac{2q(q-3)}{q-1}\lambda_1^2x_n+\left(\frac{2q(q-3)}{q-1}\lambda_1^2-6q\lambda_2^2\right)z_n-\frac{4q^2}{q-1}\lambda_1^3\left(u_n-\frac{x_n}{2q}\right)
-4q^2\frac{3\lambda_1+(q-3)\lambda_2}{q-1}\lambda_1^2\left(w_n-\frac{x_n}{2q}\right);
\end{eqnarray*}

\item when $q+1\leq i\leq 2q$,
\begin{eqnarray*}
\mathbf{E}Z_i^2=1+d\Pi_3+\frac{d(d-1)}{2}\Pi_3^2+O_q(x_n^3)
\end{eqnarray*}
where
\begin{eqnarray*}
\Pi_3&=&\mathbf{E}\left[1+2q\lambda_1\left(Y_{ij}-\frac{1}{2q}\right)+2(\lambda_2-\lambda_1)\sum_{q+1\leq
\ell\leq 2q}\left(Y_{\ell j}-\frac{1}{2q}\right)\right]^2-1
\\
&=&2q\lambda_1^2x_n+2q(\lambda_1^2+\lambda_2^2)z_n+4q^2\lambda_1^2\lambda_2\left(w_n-\frac{x_n}{2q}\right);
\end{eqnarray*}

\item for $2\leq i_1<i_2\leq q$,
\begin{equation*}
\mathbf{E}Z_{i_1}Z_{i_2}=\mathbf{E}Z_2Z_q
=1+d\Pi_4+\frac{d(d-1)}{2}\Pi_4^2+O_q(x_n^3)
\end{equation*}
where
\begin{eqnarray*}
\Pi_4&=&\mathbf{E}\left[1+2q\lambda_1\left(Y_{2j}-\frac{1}{2q}\right)+2(\lambda_1-\lambda_2)\sum_{q+1\leq\ell\leq
2q}\left(Y_{\ell j}-\frac{1}{2q}\right)\right]
\\
&&\times\left[1+2q\lambda_1\left(Y_{qj}-\frac{1}{2q}\right)+2(\lambda_1-\lambda_2)\sum_{q+1\leq\ell\leq
2q}\left(Y_{\ell j}-\frac{1}{2q}\right)\right]-1
\\
&=&-\frac{6q\lambda_1^2}{q-1}x_n-\left(\frac{6q\lambda_1^2}{q-1}+6q\lambda_2^2\right)z_n+\frac{8q^2\lambda_1^3}{(q-1)(q-2)}\left(u_n-\frac{x_n}{2q}\right)
+4q^2\frac{6\lambda_1+(3q-6)\lambda_2}{(q-1)(q-2)}\lambda_1^2\left(w_n-\frac{x_n}{2q}\right);
\end{eqnarray*}

\item for $q+1\leq i\leq 2q$,
\begin{equation*}
\mathbf{E}Z_2Z_i=\mathbf{E}Z_2Z_{q+1}
=1+d\Pi_5+\frac{d(d-1)}{2}\Pi_5^2+O_q(x_n^3)
\end{equation*}
where
\begin{eqnarray*}
\Pi_5&=&\mathbf{E}\left[1+2q\lambda_1\left(Y_{2j}-\frac{1}{2q}\right)+2(\lambda_1-\lambda_2)\sum_{q+1\leq\ell\leq
2q}\left(Y_{\ell j}-\frac{1}{2q}\right)\right]
\\
&&\times\left[1+2q\lambda_1\left(Y_{(q+1)j}-\frac{1}{2q}\right)+2(\lambda_1-\lambda_2)\sum_{1\leq\ell\leq
q}\left(Y_{\ell j}-\frac{1}{2q}\right)\right]-1
\\
&=&-\frac{2q\lambda_1^2}{q-1}x_n+\left(-\frac{2q\lambda_1^2}{q-1}+2q\lambda_2^2\right)z_n-\frac{4q^2}{q-1}\lambda_1^2\lambda_2\left(w_n-\frac{x_n}{2q}\right).
\end{eqnarray*}
\end{enumerate}

\subsection{Main Expansion of $x_{n+1}$ and $z_{n+1}$}
\indent In this section, we wish to figure out the second order
recursion relation of $x_{n+1}$ and $z_{n+1}$ by virtue of the
following identity
\begin{equation}
\label{identity}
\frac{a}{s+r}=\frac{a}{s}-\frac{ar}{s^2}+\frac{r^2}{s^2}\frac{a}{s+r}.
\end{equation}
Specifically, taking $a=Z_1$, $s=2q$ and $r=\sum_{i=1}^{2q}Z_i-2q$
in~\eqref{identity} yields
\begin{equation}
\label{xexpansion}
x_{n+1}+\frac{1}{2q}=\mathbf{E}\frac{Z_1}{\sum_{i=1}^{2q}Z_i}
=\mathbf{E}\frac{Z_1}{2q}-\mathbf{E}\frac{Z_1(\sum_{i=1}^{2q}Z_i-2q)}{(2q)^2}+\mathbf{E}\frac{Z_1}{\sum_{i=1}^{2q}Z_i}\frac{(\sum_{i=1}^{2q}Z_i-2q)^2}{(2q)^2};
\end{equation}
\begin{equation}
\label{zexpansion}
z_{n+1}+\frac{1}{2q}=\mathbf{E}\frac{Z_{q+1}}{\sum_{i=1}^{2q}Z_i}
=\mathbf{E}\frac{Z_{q+1}}{2q}-\mathbf{E}\frac{Z_{q+1}(\sum_{i=1}^{2q}Z_i-2q)}{(2q)^2}+\mathbf{E}\frac{Z_{q+1}}{\sum_{i=1}^{2q}Z_i}\frac{(\sum_{i=1}^{2q}Z_i-2q)^2}{(2q)^2}.
\end{equation}

Finally we plug the results of Section~\ref{Zresults}
into~\eqref{xexpansion} and~\eqref{zexpansion} and take
substitutions of $\mathcal{X}_n=x_n+z_n$ and $\mathcal{Z}_n=-z_n$. Thus,
there is a two dimensional recursive formula of the linear diagonal
canonical form:
\begin{equation}
\label{X}
\left\{\begin{array}{ll}
\mathcal{X}_{n+1}=d\lambda_1^2\mathcal{X}_n+\frac{d(d-1)}{2}\left(\frac{2q(q-3)}{q-1}\lambda_1^4\mathcal{X}_n^2+4q\lambda_1^2\lambda_2^2\mathcal{X}_n\mathcal{Z}_n\right)
+R_x+R_z+V_x
\\
\
\\
\mathcal{Z}_{n+1}=d\lambda_2^2\mathcal{Z}_n+\frac{d(d-1)}{2}\left(\frac{q}{q-1}\lambda_1^4\mathcal{X}_n^2-4q\lambda_2^4\mathcal{Z}_n^2\right)
-R_z+V_z
\end{array}
\right.
\end{equation}
where
\begin{equation*}
R_x=\mathbf{E}\left(\frac{Z_1}{\sum_{i=1}^{2q}Z_i}-\frac{1}{2q}\right)\frac{(\sum_{i=1}^{2q}Z_i-2q)^2}{(2q)^2},\quad
R_z=\mathbf{E}\left(\frac{Z_{q+1}}{\sum_{i=1}^{2q}Z_i}-\frac{1}{2q}\right)\frac{(\sum_{i=1}^{2q}Z_i-2q)^2}{(2q)^2},
\end{equation*}
and
$$
|V_x|, |V_z|\leq
C_Vx_n^2\left(\left|\frac{u_n}{x_n}-\frac{1}{2q}\right|+\left|\frac{w_n}{x_n}-\frac{1}{2q}\right|+x_n\right)
$$
with $C_V$ a constant depending on $q$ only.

\section{Proof of the Main Theorem}
\indent If the reconstruction problem is solvable, then $\sigma(n)$ contains
significant information on the root variable. This may be expressed
in several equivalent ways (\cite{MO1}, Proposition 14).
\begin{lemma}
\label{equivalent} The nonreconstruction is equivalent to
$$
\lim_{n\to \infty}x_n=0.
$$
\end{lemma}

In order to study the stability of dynamical system~\eqref{X}, we
expect $R_x, R_z$ and $V_x, V_z$ to be just small perturbations, for
example, of the order $o(x^2_n)$. It is known that fixed finite
different vertices far away from the root can affect the root
little, based on which, it is possible to explore further
the concentration analysis. Analogous to the concentration analysis
in~\cite{S, LN}, we can verify that
$\frac{Z_1}{\sum_{i=1}^{2q}Z_i}$,
$\frac{Z_{q+1}}{\sum_{i=1}^{2q}Z_i}$ are both sufficiently around
$\frac{1}{2q}$, and thus are able to bound $R_x$ and $R_z$
in~\eqref{X}.
\begin{lemma}
\label{R} Assume $\min\{|\lambda_1|, |\lambda_2|\}>\varrho$ for some
$\varrho>0$. For any $\varepsilon>0$, there exist $N=N(q,
\varepsilon)$ and $\delta=\delta(q, \varepsilon, \varrho)>0$ such
that if $n\geq N$ and $x_n\leq\delta$ then
\begin{equation*}
|R_x|\leq\varepsilon x_n^2\quad \textup{and}\quad
|R_z|\leq\varepsilon x_n^2.
\end{equation*}
\end{lemma}
\begin{proofsect}{Proof}
For any $\eta>0$ and $1\leq i\leq 2q$, applying Cauchy-Schwartz
inequality gives
\begin{equation*}
\begin{aligned}
&\left|\mathbf{E}\frac{Z_1}{\sum_{i=1}^{2q}Z_i}\frac{(\sum_{i=1}^{2q}Z_i-2q)^2}{(2q)^2}-\mathbf{E}\frac{1}{2q}\frac{(\sum_{i=1}^{2q}Z_i-2q)^2}{(2q)^2}\right|
\\
&\leq
\mathbf{E}\frac{(\sum_{i=1}^{2q}Z_i-2q)^2}{(2q)^2}\left|\frac{Z_1}{\sum_{i=1}^{2q}Z_i}-\frac{1}{2q}\right|
\\
&\leq\eta \mathbf{E}\left(\frac{(\sum_{i=1}^{2q}Z_i-2q)^2}{(2q)^2};
\left|\frac{Z_1}{\sum_{i=1}^{2q}Z_i}-\frac{1}{2q}\right|\leq
\eta\right)
+\mathbf{E}\left(\frac{(\sum_{i=1}^{2q}Z_i-2q)^2}{(2q)^2}\mathbf{I}\left(\left|\frac{Z_1}{\sum_{i=1}^{2q}Z_i}-\frac{1}{2q}\right|>\eta\right)\right)
\\
&\leq\eta\mathbf{E}\left(\frac{(\sum_{i=1}^{2q}Z_i-2q)^2}{(2q)^2}\right)+\left(\mathbf{E}\frac{(\sum_{i=1}^{2q}Z_i-2q)^4}
{(2q)^4}\right)^{1/2}\left(\mathbf{P}\left(\left|\frac{Z_1}{\sum_{i=1}^{2q}Z_i}-\frac{1}{2q}\right|>\eta\right)\right)^{1/2}.
\end{aligned}
\end{equation*}
We can derive from the calculation for distributional recursion that
$$\mathbf{E}\left(\frac{(\sum_{i=1}^{2q}Z_i-2q)^2}{(2q)^2}\right)\leq
C_1(q)x_n^2\quad\textup{and}\quad\mathbf{E}\frac{(\sum_{i=1}^{2q}Z_i-2q)^4}{(2q)^4}\leq
C_2(q).$$ Similar to Lemma 2.11 of ~\cite{S} and Lemma 4.3 of
~\cite{LN}, there exist $C_3=C_3(q, \eta, \varrho)$ and $N=N(q,
\eta)$ such that when $n\geq N$,
$$
\mathbf{P}\left(\left|\frac{Z_1}{\sum_{i=1}^{2q}Z_i}-\frac{1}{2q}\right|>\eta\right)\leq
C_3x_n^6.
$$
and thus there is a $C_4=C_4(q, \eta, \varrho)$ such that
\begin{eqnarray*}
|R_x|=\left|\mathbf{E}\left(\frac{Z_1}{\sum_{i=1}^{2q}Z_i}-\frac{1}{2q}\right)\frac{(\sum_{i=1}^{2q}Z_i-2q)^2}{(2q)^2}\right|
\leq\eta C_1x_n^2.
\end{eqnarray*}
Finally it suffices to take $C_1\eta=\varepsilon/2$, so if
$x_n\leq\delta$, then $R_x\leq \varepsilon x_n^2$.
Similar analysis gives $R_z\leq \varepsilon x_n^2$.
\end{proofsect}

Before investigating the concentration about $V_x$ and $V_z$, we need to firstly prove the following two lemmas. One shows that the value of $x_n$ does not drop too fast to be non-reconstruction, and the other improves the result of Lemma \ref{thm:xn} by
verifying the strict positivity of the sum of $x_n$ and  $z_n$.
\begin{lemma}
	\label{ndtf}
	For any $\varrho>0$, there exists a constant $\gamma=\gamma(q,
	\varrho)>0$ such that
	$$
	x_{n+1}\geq \gamma x_n,
	$$
	for all $n$, if $\min\{|\lambda_1|,
	|\lambda_2|\}>\varrho$,
\end{lemma}
\begin{proofsect}{Proof}
	Similarly as~\eqref{posterior}, for a configuration $A$ on $\mathbb{T}_{u_1}\cap L(n+1)$ define the posterior function
	\begin{eqnarray*}
		g_{n+1}(1, A)&=&\mathbf{P}(\sigma_\rho=1\mid\sigma_1(n+1)=A)
		\\
		&=&\frac{1}{2q}+p_0\left(f_n(1,
		A)-\frac{1}{2q}\right)+p_1\sum_{i=2}^q\left(f_n(i,
		A)-\frac{1}{2q}\right)+p_2\sum_{i=q+1}^{2q}\left(f_n(i, A)-\frac{1}{2q}\right)
		\\
		&=&\frac{1}{2q}+\lambda_1\left(f_n(1,
		A)-\frac{1}{2q}\right)+\frac{\lambda_1-\lambda_2}{q}\sum_{i=q+1}^{2q}\left(f_n(i,
		A)-\frac{1}{2q}\right)
	\end{eqnarray*}
	and then
	\begin{eqnarray*}
		\mathbf{E}g_{n+1}(1, \sigma_1^1(n+1))
		&=&\frac{1}{2q}+\lambda_1\mathbf{E}\left(Y_{11}(n)-\frac{1}{2q}\right)+\frac{\lambda_1-\lambda_2}{q}q\mathbf{E}\left(Y_{(q+1)1}(n)-\frac{1}{2q}\right)
		\\
		&=&\frac{1}{2q}+\lambda_1^2x_n+(\lambda_1^2-\lambda_2^2)z_n.
	\end{eqnarray*}
	
	The estimator that chooses a state with probability $f_{n+1}(i,
	\sigma_1(n+1))$ correctly reconstructs the root with probability
	$\frac{1}{2q}+\lambda_1^2x_n+(\lambda_1^2-\lambda_2^2)z_n$. It is
	apparent that this probability must be less than the
	maximum-likelihood estimator~\cite{MM}. Therefore we can establish the inequality:
		\begin{eqnarray*}
		\frac{1}{2q}+\lambda_1^2x_n+(\lambda_1^2-\lambda_2^2)z_n&\leq&\mathbf{E}\max_{1\leq i\leq 2q}X_i(n+1)
		\\
		&\leq&\frac{1}{2q}+\left(\mathbf{E}\max_i\left(X_i(n+1)-\frac{1}{2q}\right)^2\right)^{1/2}
		\\
		&\leq&\frac{1}{2q}+\left(\mathbf{E}\sum_{i=1}^{2q}\left(X_i(n+1)-\frac{1}{2q}\right)^2\right)^{1/2}
		\\
		&=&\frac{1}{2q}+x_{n+1}^{1/2}.
	\end{eqnarray*}
	If $\lambda_1^2\geq \lambda_2^2$, then it is concluded from $x_n+z_n\geq
	0$ in Lemma~\ref{xn} that
	\begin{eqnarray*}
		\lambda_2^2x_n
		\leq\lambda_2^2x_n+(\lambda_1^2-\lambda_2^2)(x_n+z_n)
		=\lambda_1^2x_n+(\lambda_1^2-\lambda_2^2)z_n
		\leq x_{n+1}^{1/2}.
	\end{eqnarray*}
	On the other hand, if $\lambda_1^2\leq\lambda_2^2$ then
	$\lambda_1^2x_n\leq x_{n+1}^{1/2}$ because of $z_n\leq0$. To sum up,
	we always have
	$$
	\min\{\lambda_1^2, \lambda_2^2\}x_n\leq x_{n+1}^{1/2}.
	$$

	Next choose $\varepsilon=\varrho^2$. It can be concluded from~\eqref{X}, Lemma~\ref{R}, as well as the inequalities
	\begin{equation}
	\label{unxn1} \left|\frac{u_n}{x_n}-\frac{1}{2q}\right|\leq1,\quad
	\left|\frac{w_n}{x_n}-\frac{1}{2q}\right|\leq1
	\end{equation}
	that there exists a $\delta=\delta(q, \varepsilon)>0$
	when $x_n<\delta$,
	$$
	x_{n+1}\geq (d\min\{\lambda_1^2, \lambda_2^2\}-\varepsilon)x_n\geq
	(d-1)\varrho^2x_n\geq\varrho^2x_n.
	$$
	On the other hand, if $x_n\geq \delta$ then
	$x_{n+1}\geq(\min\{\lambda_1^2,
	\lambda_2^2\}x_n)^2\geq\varrho^4\delta x_n$. Finally taking
	$\gamma=\min\{\varrho^2, \varrho^4\delta\}$ completes the proof.
\end{proofsect}

\begin{lemma}
\label{non0} Assume $\lambda_1\neq 0$. For any nonnegative $n\in
\mathbb{Z}$, we always have
$$
x_n+z_n>0.
$$
\end{lemma}
\begin{proofsect}{Proof}
In Lemma~\ref{xn} we have proved that $x_n+z_n\geq 0$, so here it
suffices to exclude the equality. Now we refer to the contradiction
and assume $x_n+z_n=0$ for some $n\in \mathbb{N}$. It follows that
if $i\neq j$ are in the same configuration set, then
\begin{eqnarray*}
\mathbf{E}(X_i(n)-X_j(n))^2=2EX_i^2(n)-2EX_i(n)X_j(n)
=\frac{x_n+z_n}{q-1} = 0.
\end{eqnarray*}
Therefore $X_1(n)=X_2(n)=\cdots=X_q(n)$ and
$X_{q+1}(n)=X_{q+2}(n)=\cdots=X_{2q}(n)$ a.s., that is, for any
configuration combination $A$ on the $n$th level, we always have
$$
\mathbf{P}(\sigma_\rho=1\mid\sigma(n)=A)=\mathbf{P}(\sigma_\rho=2\mid\sigma(n)=A).
$$
Denote the leftmost vertex on the $n$th level by $v_n(1)$, and it
follows
$$
\mathbf{P}(\sigma_\rho=1\mid\sigma_{v_n(1)}=1)=\mathbf{P}(\sigma_\rho=2\mid\sigma_{v_n(1)}=1).
$$
Define the transition matrices at distance $s$ by
$$
U_s=M_{1, 1}^s, \quad V_s=M_{1, 2}^s, \quad \textup{and} \quad
W_s=M_{1, q+1}^s,
$$
and then it is convenient to figure out the iterative formulae for
them viz.
\begin{eqnarray*}
\left\{\begin{array}{ll} U_s=p_0U_{s-1}+(q-1)p_1V_{s-1}+qp_2W_{s-1}
\\
V_s=p_1U_{s-1}+[p_0+(q-2)p_1]V_{s-1}+qp_2W_{s-1}
\\
W_s=p_2U_{s-1}+(q-1)p_2V_{s-1}+[p_0+(q-1)p_1]W_{s-1}.
\end{array}
\right.
\end{eqnarray*}
To evaluate this three order recursive system, start with the
difference of the first two equation
$$
U_s-V_s=\lambda_1(U_{s-1}-V_{s-1}),
$$
and then in light of the initial conditions $U_0=1$ and $V_0=W_0=0$,
it follows that
\begin{equation}
\label{U-V} U_s-V_s=\lambda_1^s.
\end{equation}
Finally, from the reversible property of the channel, we can conclude
that
$$
\lambda_1^n=U_n-V_n=\mathbf{P}(\sigma_\rho=1\mid\sigma_{v_n(1)}=1)-\mathbf{P}(\sigma_\rho=2\mid\sigma_{v_n(1)}=1)=0,
$$
i.e., $\lambda_1=0$, a contradiction to the assumption of
$\lambda_1\neq0$.
\end{proofsect}

The following result helps estimate the terms
$u_n-\frac{x_n}{2q}$ and $w_n-\frac{x_n}{2q}$ when $x_n$ is small.
\begin{lemma}
\label{unconcentration} Assume $|\lambda_2|>\varrho$ and
$|\lambda_1|=|\lambda_2|$ or $|\lambda_1|/|\lambda_2|\geq\kappa$ for
some $\kappa>1$. For any $\varepsilon>0$, there exist $N=N(q,
\kappa, \varepsilon)$ and $\delta=\delta(q, \kappa, \varrho,
\varepsilon)>0$ such that if $n\geq N$ and $x_n\leq\delta$ then
$$
\left|\frac{u_n}{x_n}-\frac{1}{2q}\right|<\varepsilon \quad
\textup{and} \quad
\left|\frac{w_n}{x_n}-\frac{1}{2q}\right|<\varepsilon.
$$
\end{lemma}
\begin{proofsect}{Proof}
Applying the identity~\eqref{identity} again, we have
\begin{equation}
\label{u}
\begin{aligned}
u_{n+1}&=\mathbf{E}\frac{\left(Z_1-\frac{1}{2q}\sum_{i=1}^{2q}Z_i\right)^2}{\left(\sum_{i=1}^{2q}Z_i\right)^2}
\\
&=\frac{1}{4q^2}\mathbf{E}\left(Z_1-\frac{1}{2q}\sum_{i=1}^{2q}Z_i\right)^2-\frac{1}{16q^4}\mathbf{E}\left(Z_1-\frac{1}{2q}\sum_{i=1}^{2q}Z_i\right)^2\left(\left(\sum_{i=1}^{2q}Z_i\right)^2-4q^2\right)
\\
&\quad+\frac{1}{16q^4}\mathbf{E}\frac{\left(Z_1-\frac{1}{2q}\sum_{i=1}^{2q}Z_i\right)^2}{\left(\sum_{i=1}^{2q}Z_i\right)^2}\left(\left(\sum_{i=1}^{2q}Z_i\right)^2-4q^2\right)
\end{aligned}
\end{equation}
Next estimate these expectations term by term. Again we remark that all
the $O_q$ terms in the following context only depend on $q$:
\begin{equation*}
\begin{aligned}
&\mathbf{E}\left(Z_1-\frac{1}{2q}\sum_{i=1}^{2q}Z_i\right)^2
\\
&=\mathbf{E}(Z_1-1)^2-\frac{2}{2q}\mathbf{E}(Z_1-1)\left(\sum_{i=1}^{2q}Z_i-2q\right)+\frac{1}{4q^2}E\left(\sum_{i=1}^{2q}Z_i-2q\right)^2
\\
&=2dq\lambda_1^2x_n+2dq(\lambda_1^2-\lambda_2^2)z_n+4dq^2\lambda_1^3\left(u_n-\frac{x_n}{2q}\right)
+12dq^2\lambda_1^2(\lambda_1-\lambda_2)\left(w_n-\frac{x_n}{2q}\right)+O_q(x_n^2)
\end{aligned}
\end{equation*}
and similarly,
\begin{eqnarray*}
\mathbf{E}\left(Z_1-\frac{1}{2q}\sum_{i=1}^{2q}Z_i\right)^2\left(\left(\sum_{i=1}^{2q}Z_i\right)^2-4q^2\right)&=&O_q(x_n^2),
\end{eqnarray*}
as well as
\begin{eqnarray*}
\mathbf{E}\left(\left(\sum_{i=1}^{2q}Z_i\right)^2-4q^2\right)^2&=&O_q(x_n^2).
\end{eqnarray*}
Substituting these bounds into~\eqref{u} gives
\begin{equation}
\label{unxn}
u_{n+1}=\frac{x_{n+1}}{2q}+d\lambda_1^3\left(u_n-\frac{x_n}{2q}\right)+3d\lambda_1^2(\lambda_1-\lambda_2)\left(w_n-\frac{x_n}{2q}\right)+O_q(x_n^2)
\end{equation}
by means of
$x_{n+1}=d\lambda_1^2x_n+d(\lambda_1^2-\lambda_2^2)z_n+O_q(x_n^2)$. Moreover, the similar expansion of $w_{n+1}$ would be
\begin{eqnarray*}
w_{n+1}=\frac{1}{4q^2}\mathbf{E}(Z_{q+1}-1)^2+O_q(x_n^2)
=\frac{x_{n+1}}{2q}+d\lambda_1^2\lambda_2\left(w_n-\frac{x_n}{2q}\right)+O_q(x_n^2),
\end{eqnarray*}
and thus
\begin{equation}
\label{w}
\frac{w_{n+1}}{x_{n+1}}-\frac{1}{2q}=d\lambda_1^2\lambda_2\frac{x_n}{x_{n+1}}\left(\frac{w_n}{x_n}-\frac{1}{2q}\right)+O_q\left(\frac{x_n^2}{x_{n+1}}\right).
\end{equation}

Next we consider this discussion in the $\mathcal{X}O\mathcal{Z}$ plane. First consider the
case of $\kappa>1$. Then in a small neighborhood of $(0, 0)$,
because of $d\lambda_2^2<\kappa^2d|\lambda_2^2|\leq d\lambda_1^2<1$
and $\mathcal{X}_n>0$, the discrete trajectories approach to the origin point
"tangential" to the $\mathcal{X}$-axis if $x_n$ is small enough for some
$n$~\cite{BER}. Besides, the conclusion of Lemma~\ref{non0} excludes
the trajectory along $\mathcal{Z}$-axis. Then for some $M>1$, there exist
absolute constants $N_1=N_1(q, \kappa, M)$ and $\delta_1=\delta_1(q,
\kappa, M)$ such that if $n\geq N_1$ and $x_n\leq\delta_1$, we have
simultaneously $\mathcal{X}_n\geq M\mathcal{Z}_n$ and
$$\frac{1}{M(M+1)}d\lambda_1^2x_n+O_q(x_n^2)>0,$$ where the remainder
term $O_q(x_n^2)$ comes from the expansion of $x_{n+1}$.
Consequently it follows $x_n+z_n=\mathcal{X}_n\geq \frac{M}{M+1}x_n$, which
yields, in connection with $z_n\leq 0$ in Lemma~\ref{xn},
\begin{equation}
\label{ratio}
\begin{aligned}
\frac{x_n}{x_{n+1}}=\frac{x_n}{d\lambda_1^2x_n+d(\lambda_1^2-\lambda_2^2)z_n+O_q(x_n^2)}
\leq\frac{x_n}{\frac{M}{M+1}d\lambda_1^2x_n+O_q(x_n^2)}
\leq\frac{x_n}{\left(1-\frac1M\right)d\lambda_1^2x_n}
=\frac{M}{M-1}\frac{1}{d\lambda_1^2}.
\end{aligned}
\end{equation}
The second case to be taken into account, is $|\lambda_1|=|\lambda_2|$. In
view of $1/2\leq d\lambda^2=d\lambda_1^2\leq1$, there exist also
absolute constants $N_2=N_2(q, M)$ and $\delta_2=\delta_2(q, M)$
such that if $n\geq N_2$ and $x_n\leq\delta_2$ then
\begin{eqnarray*}
\frac{x_n}{x_{n+1}}=\frac{x_n}{d\lambda_1^2x_n+O_q(x_n^2)}
\leq\frac{x_n}{\left(1-\frac1M\right)d\lambda_1^2x_n}
=\frac{M}{M-1}\frac{1}{d\lambda_1^2}.
\end{eqnarray*}

For fixed $k$, it is known from~\eqref{X} that
\begin{eqnarray*}
|x_{n+1}-(d\lambda_1^2\mathcal{X}_n+d\lambda_2^2\mathcal{Z}_n)|\leq C(q)x_n^2,
\end{eqnarray*}
and then there exists a $\delta_3=\delta_3(q, \kappa, M,
k)<\min\{\delta_1, \delta_2\}$ such that if $x_n<\delta_3$ then
$x_{n+\ell}<2\delta_3$ for any $1\leq\ell\leq k$. Therefore for any
positive integer $k$, iterating $k$ times~\eqref{w} yields
\begin{eqnarray*}
\frac{w_{n+k}}{x_{n+k}}-\frac{1}{2q}&=&d\lambda_1^2\lambda_2\frac{x_{n+k-1}}{x_{n+k}}\left(\frac{w_{n+k-1}}{x_{n+k-1}}-\frac{1}{2q}\right)+O_q\left(x_{n+k-1}\frac{x_{n+k-1}}{x_{n+k}}\right)
\\
&=&(d\lambda_1^2\lambda_2)^k\left(\prod_{\ell=1}^k\frac{x_{n+\ell-1}}{x_{n+\ell}}\right)\left(\frac{w_n}{x_n}-\frac{1}{2q}\right)+R,
\end{eqnarray*}
where
\begin{eqnarray*}
(d\lambda_1^2\lambda_2)^k\left(\prod_{\ell=1}^k\frac{x_{n+\ell-1}}{x_{n+\ell}}\right)
\leq(d\lambda_1^2|\lambda_2|)^k\left(\frac{M}{M-1}\frac{1}{d\lambda_1^2}\right)^k=
\left(\frac{M}{M-1}|\lambda_2|\right)^k
\end{eqnarray*}
and
$$
|R|\leq
2C\delta_3\left(\sum_{i=1}^k\left(\frac{M}{M-1}\frac{1}{d\lambda_1^2}\right)^i(d\lambda_1^2|\lambda_2|)^{i-1}\right)
\leq\delta_3\frac{1-\left(\frac{M}{M-1}|\lambda_2|\right)^k}{1-\left(\frac{M}{M-1}|\lambda_2|\right)}\frac{M}{M-1}\frac{1}{d\lambda_1^2}
$$
with $C$ denoting the $O_q$ constant in~\eqref{w}. From the identity
(i) in Lemma~\ref{xnun}, it is easy to obtain
$0\leq\frac{w_n}{x_n}\leq\frac{1}{q}$, which implies
$$\left|\frac{w_n}{x_n}-\frac{1}{2q}\right|\leq \frac{1}{2q}.$$ Noticing the fact of $|\lambda_2|\leq |\lambda_1|\leq d^{-1/2}\leq
1/\sqrt{2}$, it is possible to achieve $\frac{M}{M-1}|\lambda_2|<1$
by choosing arbitrary $M\geq 4$, say, $M=4$. Therefore it is
feasible to take $k=k(\varepsilon)$ sufficiently large and
$\delta_4=\delta_4(q, \kappa, k, \varepsilon)=\delta_4(q, \kappa,
\varepsilon)<\delta_3$ sufficiently small to guarantee
$$\left|\frac{w_{n+k}}{x_{n+k}}-\frac{1}{2q}\right|<\varepsilon.$$
Finally in view of $|\lambda_2|>\varrho$, there exists
$\gamma=\gamma(q, \varrho)$ by Lemma~\ref{ndtf} satisfying
$x_{n-k}\leq \gamma^{-k}x_n$, and thus by choosing $N=N(q, \kappa,
\varepsilon, k)=N(q, \kappa, \varepsilon)>\max\{N_1+k, N_2+k\}$ and
$\delta=\gamma^k\delta_4$, if $x_n\leq\delta$ and $n\geq N$ then
\begin{equation}
\label{wnxn} \left|\frac{w_n}{x_n}-\frac{1}{2q}\right|< \varepsilon.
\end{equation}
Finally the second part of the lemma follows by
plugging~\eqref{wnxn} into~\eqref{unxn} and proceeding similarly as
above.
\end{proofsect}

\begin{proofsect}{Proof of the Main Theorem}
First for any fixed $\varrho>0$, consider
$\varrho<|\lambda_2|<|\lambda_1|$. By Lemma~\ref{equivalent}, it
suffices to establish that when $d\lambda_1^2$ is close enough to
$1$, $\mathcal{X}_n$ does not converge to $0$. Because it implies that $x_n$ does not
converge to $0$ either, considering $0\leq \mathcal{X}_n=x_n+z_n\leq x_n$.
Therefore it is convenient to make $|\lambda_2|>\varrho$ fixed and
just $\lambda_1$ varying, and then without loss of generality,
assume $d\lambda_1^2>\frac{1+d\lambda_2^2}{2}$. Consequently choose
$\kappa=\kappa(d,
\lambda_2)=\left(\frac{1+d\lambda_2^2}{2d\lambda_2^2}\right)^{1/2}>1$
and thus $|\lambda_1|/|\lambda_2|\geq \kappa$.

As in Lemma~\ref{unconcentration}, display our proof
in the $\mathcal{X}O\mathcal{Z}$ plane. With the condition of $q\geq 4$ and~\eqref{X},
it is apparent that
\begin{eqnarray*}
\mathcal{X}_{n+1} &=&
d\lambda_1^2\mathcal{X}_n+\frac{d(d-1)}{2}\left(\frac{2q(q-3)}{q-1}\lambda_1^4\mathcal{X}_n^2+4q\lambda_1^2\lambda_2^2\mathcal{X}_n\mathcal{Z}_n\right)
+R_x+R_z+V_x
\\
&\geq&d\lambda_1^2\mathcal{X}_n+\frac{d(d-1)}{2}\frac{2q(q-3)}{q-1}\lambda_1^4\mathcal{X}_n^2-|R_x|-|R_z|
-C_Vx_n^2\left(\left|\frac{u_n}{x_n}-\frac{1}{2q}\right|+\left|\frac{w_n}{x_n}-\frac{1}{2q}\right|+x_n\right),
\end{eqnarray*}
where the last inequality comes from $|\lambda_1|\leq d^{-1/2}<1$.
Then by Lemma~\ref{unconcentration} and Lemma~\ref{R}, there exist
$N=N(q, \kappa, \varrho)$ and $\delta=\delta(q, d, \kappa,
\varrho)>0$ such that if $n\geq N$ and $x_n\leq\delta$, then in the
small neighborhood of the origin point $(0, 0)$, we have
$\mathcal{X}_n\geq \mathcal{Z}_n$ and thus
$\mathcal{\mathcal{X}}_n\geq\frac{x_n}2$. Meanwhile, the following estimates hold simultaneously:
$$
x_n\leq
\frac{1}{48C_V}\frac{d(d-1)}{2}\frac{2q(q-3)}{q-1}\lambda_1^4;
$$
$$
\left|\frac{u_n}{x_n}-\frac{1}{2q}\right|, \left|\frac{w_n}{x_n}-\frac{1}{2q}\right|\leq
\frac{1}{48C_V}\frac{d(d-1)}{2}\frac{2q(q-3)}{q-1}\lambda_1^4;
$$
$$
|R_x|, |R_z|\leq
\frac{1}{32}\frac{d(d-1)}{2}\frac{2q(q-3)}{q-1}\lambda_1^4x_n^2\leq
\frac{1}{8}\frac{d(d-1)}{2}\frac{2q(q-3)}{q-1}\lambda_1^4\mathcal{X}_n^2.
$$
Therefore, the quadratic term of $\mathcal{X}_n^2$ is big enough to control the remainder terms:
\begin{equation}
\label{recursiveinequality}
\begin{aligned}
\mathcal{X}_{n+1}\geq
d\lambda_1^2\mathcal{X}_n+\frac12\frac{d(d-1)}{2}\frac{2q(q-3)}{q-1}\lambda_1^4\mathcal{X}_n^2
=\mathcal{X}_n\left[d\lambda_1^2+\frac12\frac{d(d-1)}{2}\frac{2q(q-3)}{q-1}\lambda_1^4\mathcal{X}_n\right].
\end{aligned}
\end{equation}

Take $\varepsilon=\min\{\frac14\gamma^N, \gamma\delta\}> 0$, where
$\gamma=\gamma(q, \varrho)>0$ is the constant in Lemma~\ref{ndtf}. Because $q\geq 4$ and $\varepsilon$ is independent of
$\lambda_1$, we can choose $|\lambda_1|<d^{-1/2}$ to make
\begin{equation}
\label{inequality}
d\lambda_1^2+\frac12\frac{d(d-1)}{2}\frac{2q(q-3)}{q-1}\lambda_1^4\varepsilon>1.
\end{equation}

Since $x_0=1-\frac{1}{2q}>\frac12$, it is concluded that
$x_n>\frac12\gamma^n\geq 2\varepsilon$ when $n\leq N$, in addition,
$\mathcal{X}_N\geq
\frac{\mathcal{X}_N+\mathcal{Z}_N}{2}=\frac{x_N}{2}\geq
\varepsilon$. Now suppose $\mathcal{X}_n\geq \varepsilon$ for some
$n\geq N$. Then display our discussion of $\mathcal{X}_n$ as
follows:
\begin{enumerate}[(1)]
\item If $\mathcal{X}_n\geq 2\gamma^{-1}\varepsilon$, then
$$
\mathcal{X}_{n+1}\geq\frac{x_{n+1}}{2}\geq\frac{\gamma
x_n}{2}\geq\frac{\gamma \mathcal{X}_n}{2}\geq\varepsilon;
$$
\item If $\varepsilon\leq \mathcal{X}_n\leq 2\gamma^{-1}\varepsilon$, then
$x_n\leq \frac{\mathcal{X}_n}{2}\leq\gamma^{-1}\varepsilon\leq
\delta$, and thus it follows by~\eqref{recursiveinequality}
and~\eqref{inequality} that
\begin{eqnarray*}
x_{n+1}\geq \mathcal{X}_{n+1} \geq
\mathcal{X}_n\left[d\lambda_1^2+\frac12\frac{d(d-1)}{2}\frac{2q(q-3)}{q-1}\lambda_1^4\mathcal{X}_n\right]
\geq \mathcal{X}_n \geq\varepsilon.
\end{eqnarray*}
\end{enumerate}
Finally show by induction that for all $n$ that $x_n\geq
\mathcal{X}_n\geq\varepsilon$. Consequently it is established that
the Kesten-Stigum bound is not tight.

The second case to be considered is $|\lambda_1|=|\lambda_2|$,
under which there are two equal multipliers in this nonlinear second
order point mapping, the origin point must be a star node. Although
the principal axis is undetermined, just by the comparison of the
quadratic terms and $q\geq4$, it is concluded that
\begin{equation*}
\begin{aligned}
&\frac{d(d-1)}{2}\left(\frac{2q(q-3)}{q-1}\lambda_1^4\mathcal{X}_n^2+4q\lambda_1^4\mathcal{X}_n\mathcal{Z}_n\right)-\frac{d(d-1)}{2}\left(\frac{q}{q-1}\lambda_1^4\mathcal{X}_n^2-4q\lambda_1^4\mathcal{Z}_n^2\right)
\\
&=\frac{d(d-1)}{2}\left(\frac{2q^2-7q}{q-1}\lambda_1^4\mathcal{X}_n^2+4q\lambda_1^4\mathcal{X}_n\mathcal{Z}_n+4q\lambda_1^4\mathcal{Z}_n^2\right)
\\
&\geq\frac{d(d-1)}{2}\lambda_1^4x_n^2,
\end{aligned}
\end{equation*}
and thus the decay rate of $\mathcal{X}_n$ is much slower than
$\mathcal{Z}_n$ if $x_n$ is sufficiently small. Therefore in light
of the preceding discussion, there still exist $N=N(q)$ and
$\delta=\delta(q)$ such that if $n\geq N$ and $x_n\leq\delta$, we
have $\mathcal{X}_n\geq \mathcal{Z}_n$ and thus
$x_n=\mathcal{X}_n+\mathcal{Z}_n\leq2\mathcal{X}_n$. Then the rest
would be the same as the first part.
\end{proofsect}

\section*{Acknowledgements}
We would like to thank two anonymous reviewers who
provided us with many constructive and helpful comments.

\end{document}